\documentclass[12pt]{article}

\usepackage{amsmath,amsthm,amsfonts,amssymb,amscd,latexsym}
\usepackage{float}
\usepackage{graphicx}
\begin{document}

\theoremstyle{plain} 
\newtheorem{thm}{\sc Theorem}
\newtheorem{asser}{\sc Assertion}
\newtheorem{lem}[thm]{\sc Lemma}
\newtheorem{cor}{\sc Corollary}
\newtheorem{prop}[thm]{\sc Proposition}

\theoremstyle{remark}
\newtheorem{Case}{\bf Case}
\newtheorem{rem}[thm]{\sc Remark}

\theoremstyle{definition}
\newtheorem{Def}{\bf Definition}
\newtheorem*{pf}{\sc Proof}

\theoremstyle{definition}
\newtheorem{Conj}{\bf Conjecture}

\newcommand{\N}{{\mathbb N}}
\newcommand{\ve}{{\varepsilon}}
\newcommand{\vk}{{\varkappa}}
\newcommand{\be}{{\beta}}
\newcommand{\la}{{\lambda}}
\newcommand{\teta}{{\theta}}
\newcommand{\La}{{\Lambda}}
\newcommand{\ga}{{\gamma}}
\newcommand{\Ga}{{\Gamma}}
\newcommand{\po}{{\partial}}
\newcommand{\ov}{\overline}
\newcommand{\re}{{\mathbb{R}}}
\newcommand{\Pc}{{\mathbb{P}}}
\newcommand{\dc}{{\mathcal D}}
\newcommand{\na}{{\mathbb N}}
\newcommand{\co}{{\mathbb C}}
\newcommand{\fc}{{\mathcal F}}
\newcommand{\al}{{\alpha}}
\newcommand{\lc}{{\mathcal L}}
\newcommand{\om}{{\omega}}
\newcommand{\OM}{{\bf \Omega}}
\newcommand{\Z}{{\mathbb Z}}
\newcommand{\rth}{{\mathbb{R}^3}}
\newcommand{\W}{{\mathcal W}}
\renewcommand{\H}{{\mathcal H}}
\newcommand{\cc}{{\mathcal C}}
\newcommand{\E}{{\mathcal E}}
\newcommand{\PP}{{\mathcal{P}}}
\newcommand{\M}{{\mathcal M}}
\newcommand{\n}{{\noindent}}
\newcommand{\m}{{\underline{m}}}
\newcommand{\UO}{{\underline{O}}}
\newcommand{\uo}{{\underline{o}}}
\newcommand{\A} {{\bf A}}


%
\def\theequation{\arabic{equation}}


\DeclareGraphicsExtensions{.ps,.jpg}


\def\bba{{\Bbb A}}
\def\bbb{{\Bbb B}}
\def\bbc{{\Bbb C}}
\def\bbd{{\Bbb D}}
\def\bbe{{\Bbb E}}
\def\bbf{{\Bbb F}}
\def\bbg{{\Bbb G}}
\def\bbh{{\Bbb H}}
\def\bbi{{\Bbb I}}
\def\bbj{{\Bbb J}}
\def\bbk{{\Bbb K}}
\def\bbl{{\Bbb L}}
\def\bbm{{\Bbb M}}
\def\bbn{{\Bbb N}}
\def\bbo{{\Bbb O}}
\def\bbp{{\Bbb P}}
\def\bbq{{\Bbb Q}}
\def\bbr{{\Bbb R}}
\def\bbs{{\Bbb S}}
\def\bbt{{\Bbb T}}
\def\bbu{{\Bbb U}}
\def\bbv{{\Bbb V}}
\def\bbw{{\Bbb W}}
\def\bbx{{\Bbb X}}
\def\bby{{\Bbb Y}}
\def\bbz{{\Bbb Z}}

\def\SC{\mathcal C}
\def\SE{\mathcal E}
\def\SF{\mathcal F}
\def\SM{\mathcal M}
\def\SO{\mathcal O}
\def\SP{\mathcal P}
\def\SR{\mathcal R}
\def\ST{\mathcal T}
\def\SU{\mathcal U}
\def\SH{\mathcal H}
\def\SY{\mathcal Y}

\def\a{\alpha}
\def\b{\beta}

\def\cd{\cdot}
\def\e{\epsilon}
\def\ve{\varepsilon}
\def\equdef{\overset{\text{def}}{=}}
\def\f{\frac}
\def\g{\gamma}
\def\G{\Gamma}
\def\hra{\hookrightarrow}
\def\k{\kappa}
\def\lb\{{\left\{}
\def\la{\lambda}
\def\La{\Lambda}
\def\lla{\longleftarrow}
\def\lm{\limits}
\def\lra{\longrightarrow}
\def\dllra{\Longleftrightarrow}
\def\llra{\longleftrightarrow}
\def\n{\nabla}
\def\ngth{\negthickspace}
\def\ola{\overleftarrow}
\def\Om{\Omega}
\def\om{\omega}
\def\op{\oplus}
\def\oper{\operatorname}
\def\oplm{\operatornamewithlimits}
\def\ora{\overrightarrow}

\def\ova{\overarrow}
\def\ox{\otimes}
\def\p{\partial}
\def\rb\}{\right\}}
\def\s{\sigma}
\def\sbq{\subseteq}
\def\spq{\supseteq}
\def\sqp{\sqsupset}
\def\supth{{\text{th}}}
\def\T{\Theta}
\def\th{\theta}
\def\tl{\tilde}
\def\thra{\twoheadrightarrow}
\def\un{\underline}
\def\ups{\upsilon}
\def\vp{\varphi}
\def\wh{\widehat}
\def\wt{\widetilde}
\def\x{\times}
\def\z{\zeta}
\def\({\left(}
\def\){\right)}
\def\[{\left[}
\def\]{\right]}
\def\<{\left<}
\def\>{\right>}
\def\x{\times}

\def\ga{Gau\ss\ }
\def\wei{Weierstra\ss\ }

\def\height{\mathcal H}
\def\lb\{{\left\{}
\def\rb\}{\right\}}
\def\ovt{\overline{\mathcal T}}
\def\BC{\ensuremath{\mathbf C}}
\def\BE{\ensuremath{\mathbf E}}
\newcommand{\ms}{minimal surface}
\newcommand{\mss}{minimal surfaces}
\newcommand{\BR}{\ensuremath{\mathbf R}}
\newcommand{\area}{\operatorname{Area}}
\newcommand{\ext}{\operatorname{Ext}}
\newcommand{\genus}{\operatorname{genus}}
\newcommand{\ree}{\operatorname{Re}}
\newcommand{\dist}{\operatorname{dist}}

\newcommand{\halo}[1]{\Int(#1)}
\def\Int{\mathop{\rm Int}}
\def\Re{\mathop{\rm Re}}
\def\Im{\mathop{\rm Im}}
\def\bfR{\mathop{\bf R}}
\def\bfC{\mathop{\bf C}}
\newcommand{\union}{\cup}
\newcommand{\goesto}{\rightarrow}
\def\BA{{\bold A}}
\def\BB{{\bold B}}
\def\BC{{\bold C}}
\def\BD{{\bold D}}
\def\BE{{\bold E}}
\def\BF{{\bold F}}
\def\BG{{\bold G}}
\def\BH{{\bold H}}
\def\BI{{\bold I}}
\def\BJ{{\bold J}}
\def\BK{{\bold K}}
\def\BL{{\bold L}}
\def\BM{{\bold M}}
\def\BN{{\bold N}}
\def\BO{{\bold O}}
\def\BP{{\bold P}}
\def\BQ{{\bold Q}}
\def\BR{{\bold R}}
\def\BS{{\bold S}}
\def\BT{{\bold T}}
\def\BU{{\bold U}}
\def\BV{{\bold V}}
\def\BW{{\bold W}}
\def\BX{{\bold X}}
\def\BY{{\bold Y}}
\def\BZ{{\bold Z}}

\def\per{\mathop{\rm Per}}
\def\res{\mathop{\rm Res}}
\def\Ext{\mathop{\rm Ext}}
\def\supp{\mathop{\rm supp}}
\def\id{\mathop{\rm Id}}

\def\ogup{\mathop{\Omega_{gdh}}}
\def\ogdn{\mathop{\Omega_{g^{-1}dh}}}
\newcommand{\sgn}{\operatorname{sgn}}

\newcommand{\bdy}{\partial}
\newcommand{\lap}{\bigtriangleup}
\def\({\left(}
\def\){\right)}
\def\[{\left[}
\def\]{\right]}
\def\ogup{\mathop{\Omega_{gdh}}}
\def\ogdn{\mathop{\Omega_{g^{-1}dh}}}
\def\R3/{${\bf R^3}$}
\def\Hel/{${\cal H}$}
\def\HEone/{${\cal H}e_{_1}$}
\def\Hone/{${\cal H}_{_1}$}
\def\Hk/{${\cal H}_{_k}$}
\def\He/{${\cal H}$}
\def\T/{${\cal T}_{_1}$}
\def\t/{$\cal T$}
\def\tkd/{${\cal T}_k(d)$}
\def\tec{Teichm\"uller\ }
\def\wei{Weierstrass\ }
\def\hksk{{\cal H}_{_k}/\sigma_k}
\def\s_k{\sigma_k}
\def\O/{${\cal O}$}
\def\aR/{${\cal R}$}
\newtheorem{theorem}{Theorem}
\newtheorem{assertion}{Assertion}
\newtheorem{proposition}{Proposition}
\newtheorem{remark}{Remark}
\newtheorem{lemma}{Lemma}
\newtheorem{definition}{Definition}
\newtheorem{defn}{Definition}
\newtheorem{corollary}{Corollary}
\newtheorem{observation}{Observation}
\newtheorem{conjecture}{Conjecture}
\newtheorem{question}{Question}
\newtheorem{example}{Example}

\newenvironment{proofspec}[1]%
{\smallskip\noindent{\bf Proof of #1.}\hskip \labelsep}%
{\nobreak\hfill\hfill\nobreak\copy\qedbox\par\medskip}
\newenvironment{acknowledgements}{\smallskip\noindent{\bf
Acknowledgements.}%
\hskip\labelsep}{}

\title{Minimal surfaces with the area growth of two planes; the case
  of infinite symmetry}
\author{William H. Meeks III\footnote{partially supported by NSF grant DMS-0405836.}
\\
Department of Mathematics\\University of Massachusetts\\
Amherst, MA 01003\\
\and
Michael Wolf\footnote{partially supported by NSF grants
DMS-9971563 and DMS-0139887.  Any opinions, findings, and conclusions or recommendations
 expressed in this publication are those of the authors and do not
 necessarily reflect the views of the NSF.}\\
Department of Mathematics\\   Rice University\\
Houston TX 77005
}

\maketitle

\section{Introduction.}

Consider a properly immersed minimal surface $M$ in $\rth$ with area
$A(r)$ in balls $B(r)$ of radius $r$ centered at the origin.  By the
monotonicity formula, the function $\overline{A}(r)= \frac{A(r)}{r^2}$
is monotonically increasing.  We say that $M$ has {\it area growth constant}
$A(M)\in(0, \infty]$, if $A(M)=\lim_{r \rightarrow \infty} \overline{A}(r)$.  Note that under a rigid motion or homothety $M^\prime$
of $M$, the number $A(M)=A(M^\prime)$, and so, $A(M) \geq \pi$, which is
the area growth constant of a plane.  We say
that $M$ has {\it quadratic area growth}, if $A(M)< \infty$.

Basic results in geometric measure theory imply that for any $M$ with
quadratic  area growth and for any sequence of positive numbers $t_i \rightarrow 0$,
the sequence homothetic shrinkings $M(i)= t_i M$ of $M$ contains a
subsequence that converges on compact subsets of $\rth$ to a limit
minimal cone $C$ in $\rth$ over a geodesic integral varifold in the unit
sphere $S^2$, which consists of a balanced finite
configuration of geodesic arcs with positive integer multiplicities.
$C$ is called a {\it limit tangent cone at infinity} to $M$.

In 1834, Scherk \cite{sche1} discovered a singly-periodic embedded minimal
surface $S_{\frac{\pi}{2}}$ in $\rth$ with quadratic area growth
constant $2 \pi$. Away from the $x_3$-axis, Scherk's surface is
asymptotic to the union of the $(x_1, x_3)$-plane and the $(x_2,
x_3)$-plane.  Geometrically Scherk's singly-periodic surface may be
viewed as the desingularization of these two orthogonal planes, which
form its unique limit tangent cone at infinity.  In
1988, Karcher \cite{ka4} defined a one-parameter deformation $S_\theta,\;
\theta \in (0, \frac{\pi}{2}],$ of Scherk's original surface
$S_{\frac{\pi}{2}}$, which are also called Scherk surfaces and which may be viewed as the desingularization of
two vertical planes with an angle $\theta$ between them. 
The limit tangent cone at infinity to $S_\theta$ consists of the union
of these planes.  We remark that under appropriate homothetic
scalings, the surfaces $S_\theta$ converge to a catenoid as $\theta
\rightarrow 0$.  Note that a catenoid has a plane of multiplicity two as its limit
tangent cone at infinity.

In \cite{me28}, Meeks presented the following three conjectures
related to minimal surfaces with quadratic area growth.

\begin{Conj} (Unique Limit Tangent Cone Conjecture, see Conjecture 11 in
  \cite{me28}.)  A properly immersed minimal surface in $\rth$ of quadratic area growth
has a unique limit tangent cone at infinity.
\end{Conj}

\begin{Conj} (Quadratic Area Growth Conjecture, see Conjecture 13 in
  \cite{me28}.)  A properly immersed minimal surface $M$ in $\rth$ has quadratic area
growth if and only if there exists a standard double cone in $\rth$
that intersects $M$ in a compact set.  By standard double cone, we
mean the union of all lines in $\rth$ passing through the origin and
through some $\ve$-disk
on $S^2, \;\ve < \frac{\pi}{2}$. 
\end{Conj}

\begin{Conj} (Scherk Uniqueness Conjecture, see Conjecture 10 in \cite{me28}.)
A connected properly immersed minimal surface $M$ in $\rth$ with quadratic area
growth constant $A(M) < 3 \pi$ must be a plane, a catenoid
or a Scherk singly-periodic minimal surfaces $S_\theta, \theta \in (0, \frac{\pi}{2}]$.
\end{Conj}

The main goal of this paper is to prove Conjecture 3 under the additional
hypothesis of infinite symmetry.

\begin{thm}  \label{thm1}  A connected
  properly immersed minimal surface in $\rth$ with infinite symmetry group
  and area growth constant less than $3 \pi$ is a plane, a catenoid or a Scherk
  singly-periodic minimal surface.
\end{thm}

We view Theorem \ref{thm1} as a first
step in resolving Conjecture 3. We hope that some of the ideas used in the proof of Theorem~\ref{thm1} will eventually lead to a proof
of this more general conjecture
and that such a proof will in turn lead to a positive
solution of the following fundamental singularities question.

\begin{Conj} (Isolated Singularities Conjecture, see Conjecture
  4 in \cite{me28}.)  Suppose $M$ is a minimal surface in a 
closed geodesic ball
  $B$ in a Riemannian three-manifold such that $\partial M \subset \partial B$
  and $M$ is properly embedded in $B$ punctured at the center of the
  ball.  Then, $M$ extends across the puncture to a smooth compact embedded
  minimal surface in $B$.
\end{Conj}

Suppose now that $M$ is a properly immersed minimal surface in $\rth$
with infinite symmetry group.  Then, $M$ is either a surface of revolution,
and so, is a catenoid, or $M$ is invariant under a screw motion
symmetry with possibly trivial
rotational part and the symmetry acts in an orientation preserving manner on
$M$.  In \cite{mr2}, Meeks and Rosenberg studied  properly embedded
minimal surfaces $M$ in $\rth$, which are invariant under a group
$\mathbb{Z}$ of
isometries generated by a screw motion symmetry with vertical axis 
 and whose quotient
surface $\overline{M}$ in the flat three-manifold $\rth/\mathbb{Z}$ has finite
topology.  They proved that the ends of such a $\overline{M}$ are asymptotic to horizontal
planes, vertical flat half annuli (quotients of half planes in $\rth$)
or helicoid ends in $\rth/\mathbb{Z}$.
Thus, if $M$ has quadratic area growth and $\overline{M}$ has finite
topology (finitely generated fundamental group), then $\overline{M}$ must have ends asymptotic to vertical half annuli.
Such annular ends of $\overline{M}$ are called {\it Scherk ends}, because the
 singly-periodic quotients $\overline{S}_\theta$ of the classical
Scherk examples, $S_\theta \subset \rth, \theta
\in (0, \frac{\pi}{2}],$ have this type of end.  It follows that $M$
is also invariant under a pure translation, and so, after a rigid motion and
homothety, we will assume that $\mathbb{Z}=\{(0,0,n) \mid n \in
\mathbb{Z}\}$ acts by translation on $\rth$ in this case.  

The next
theorem is the key result that we need to prove Theorem 1.

\begin{thm}  \label{thm2} The Scherk minimal surfaces 
  $\overline{S}_\theta$ are the unique connected minimal surfaces in
  $\rth/\mathbb{Z}$ with four Scherk ends.
\end{thm}

To our knowledge, this theorem is the first uniqueness result for 
singly-periodic Scherk surfaces of  genus greater than $1$.  The case
of genus $0$ was shown by Meeks and Rosenberg \cite{mr3}, and the case
of genus $1$ was carried out in the doctorial thesis of 
Hai-Ping Luo \cite{lu1}.  

We emphasize a corollary of Theorem~\ref{thm2}, by noting that 
if $\overline{M} \subset \rth/\mathbb{Z}$ has
ends asymptotic to vertical half annuli, then its lift to $\rth$,
in the complement of a vertical cylinder, is asymptotic to four
half planes. Thus, we may regard such a surface as a periodic minimal
desingularization of the intersection of two planes, and rephrase 
Theorem~\ref{thm2} as:

\begin{cor} The Scherk minimal surfaces 
  $\overline{S}_\theta$ are the unique connected periodic 
minimal desingularizations of the intersections of two planes.
\end{cor}

Our proof of Theorem \ref{thm1} is broken up into a series of 
propositions, which
appear in sections of the manuscript.  Let $\overline{M} \subset \rth/
\mathbb{Z}$ be a connected properly immersed minimal surface with four
Scherk ends. The lifted surface $M\subset \rth$ then has quadratic area
growth constant $2 \pi$.  Applying the monotonicity formula for area 
 to $M$ at a possible point
of self-intersection, we see
that $\overline{M}$ is an embedded minimal surface.
  In section 2, we prove that the
corresponding $M$ has two vertical planes of Alexandrov symmetry,
just as the classical Scherk examples have.  These planes decompose $M$ into
four nonempty closed simply connected regions; we study
two one-forms underlying the \wei representation of $M$. These forms
naturally define flat structures on each of these four regions in
$M$, and these flat structures develop to one of the two
complements
of a zigzag in $\mathbb{E}^2$.
In section 3, we prove a local rigidity theorem for the surfaces
$M$ with fixed angle $\theta$ between their Scherk ends.  The implicit
function theorem then asserts that in terms of the angle map between
 the Scherk ends $\theta \colon \mathcal{M}_{\vk}
\rightarrow (0, \frac{\pi}{2}],$ defined on the moduli space
$\mathcal{M}_{\vk}$ of  examples $\overline{M}$ of genus $\vk$ in
$\rth/\mathbb{Z}$ (defined up to congruence), we have that every
component of $\mathcal{M}_{\vk}$ is a nontrivial curve $C$ and $\theta|_C \colon C
\rightarrow \theta (C) \subset (0, \frac{\pi}{2}]$ is a
diffeomorphism.  In section 4, we prove that $\theta \colon
\mathcal{M}_{\vk} \rightarrow (0, \frac{\pi}{2}]$ is proper, and so,
$\theta|_C \colon C \rightarrow (0, \frac{\pi}{2}]$ is a
diffeomorphism.  In section 5, we prove that for $\alpha$
close to $0, \;\theta^{-1}(\alpha)$ is one of the Scherk examples
$\overline{S}_\alpha$.  Thus, $\mathcal{M}_{\vk}$ contains only one
component, which is the component of Scherk examples.  This result
proves Theorem 2.

In section 6, we prove that if $M$ is a properly embedded minimal
surface in $\rth$ with $A(M)<3 \pi,$ then $A(M)=\pi$ and $M$ is a plane
or $A(M)=2 \pi$.  Under the assumption that $M$ has infinite symmetry
group, we then prove that $M$ is either a catenoid or $M$ is invariant
under  a group $\mathbb{Z}$ of translations with $\overline{M}=M/\mathbb{Z}$ having
finite topology.  Then, our results from section 5 complete the proof
of Theorem \ref{thm1}.

Our basic strategy
of proving Theorem 2 is to show that the 
angle map $\theta$ on the moduli space is open and
proper and that examples with small angle in
$\mathcal{M}_{\vk}$ are Scherk; this strategy is 
motivated by the proofs of two previous
uniqueness theorems in the literature.
Lazard-Holly and Meeks \cite{lhm} used this approach in their
characterization of 
 the family of Scherk doubly-periodic
 minimal surfaces $\widetilde{S}_\theta,\; \theta \in 
 (0, \frac{\pi}{2}],$ which are also parametrized by the angle between
 the ends of their quotient surfaces,  as being the only properly embedded minimal surfaces in $\rth$
 with genus $0$ quotients.   A similar approach was also used by Meeks, Perez
and Ros \cite{mpr1} to characterize of the one-parameter family of
Riemann minimal examples, the helicoid and the plane, as being the only
properly embedded periodic genus $0$ minimal surfaces in $\rth$.
In another direction, Perez and Traizet \cite{PeTra1} have recently
classified the properly embedded singly-periodic minimal surfaces with quotient
surfaces having genus $0$ and finite
topology; their classification
theorem has similar structural attributes and they prove that these surfaces are precisely the Scherk
towers defined earlier by Karcher \cite{ka4}.  Their 
classification result then leads to the classification \cite{PeRoTra1} of
properly embedded doubly-periodic minimal surfaces in $\rth$ whose
quotients have genus $1$ and parallel annular ends in $\mathbb{T}
\times \re$, where $\mathbb{T}$ is a flat torus.  We remark that this last
classification result implies that these genus $1$ minimal
surfaces are the same examples which were defined by Karcher in
\cite{ka6}. 

  Finally, we remark that the idea used in our proof of
Theorem \ref{thm1} of demonstrating the local rigidity of a minimal surface in a 
moduli space of flat structures is a cornerstone of the work 
of Weber and Wolf (\cite{wolf3,ww2,ww1})
in their Teichm\"{u}ller-theoretic approach to existence
problems in minimal surface theory.

\section{Existence of Alexandrov planes of reflectional symmetry.}

The following proposition is well-known to experts in the field and
the proof (unpublished) we give is due to Harold Rosenberg.

\begin{prop} \label{symmetry} Let $M \subset \rth$ be a properly embedded
  minimal surface invariant under translation by the vectors
  $\mathbb{Z}=\{(0,0,n) \mid n \in \mathbb{Z}\}$ and such that the
  quotient surface $\overline{M}=M/\mathbb{Z}$ has four Scherk ends
  and the genus of $M$ is $\vk$.  Let $A_{x_1}$ and $A_{x_2}$ be the vertical
annuli parallel to the $x_1$ and $x_2$ axes, respectively, which are
quotients of the vertical $(x_1,x_3)$- and the $(x_2,x_3)$-planes,
respectively by the $\mathbb{Z}$-action.

Then: 
\begin{enumerate}
\item After rigid motion, we may assume $\overline{M}$ is invariant
  under reflection in the vertical annuli $A_{x_1}, A_{x_2},$ 
which each intersect
  $\overline{M}$ orthogonally in $\vk+1$ strictly convex simple closed curves contained in
  the respective annuli.
\item The regions of $\overline{M}$ on either side of $A_{x_1}$ or
  $A_{x_2}$ are graphs over their projections to the respective
  annuli.  In particular, the $(x_1,x_3)$-plane and the $(x_2,x_3)$-plane
  are Alexandrov planes of reflexive symmetry for $M$, after a
  rigid motion of $M$.  
\end{enumerate} 

\end{prop}

\proof Consider the flux vectors
$$v_i= \int_{\gamma_i}(\nabla x_1, \; \nabla x_2),\; i = 1,2,3,4$$
defined for oriented loops $\gamma_1, \gamma_2, \gamma_3, \gamma_4$ around the four
cyclically ordered Scherk ends of $\overline{M}$.  Since 
each of these vectors is a unit vector and the sum of
these vectors is zero by the divergence theorem, 
we see that after a rotation of
$M$ around the $x_3$-axis, 
we must have $v_1=(\cos \alpha, \sin \alpha), v_2=(\cos
\alpha, - \sin \alpha), v_3= -v_1$ and $v_4=-v_2$, where $\alpha \in
(0, \frac{\pi}{4}]$.  Note that the Scherk
ends of $\overline{M}$ are asymptotic to ends $A_1, A_2, A_3, A_4,$ of
flat vertical annuli, where $v_i$ is parallel to $A_i$.

Now consider the family $E(1,t)$ of vertical annuli, which are parallel
to $A_{x_1}$ and pass through
the point $(0,t,0)$, for $t \in \re$.  Let $R(1,t) \colon \re^2 \times
\re /\mathbb{Z}\rightarrow \re^2 \times \re/\mathbb{Z}$ denote reflection across $E(1,t)$, let
$\overline{M}_+(1,t)$ denote the portion of $\overline{M}$ on the side of
$E(1,t)$ which contains large positive $x_1$ valued points of
$\overline{M}$ and let $\ov{M}_-(1,t)=R(1,t)(\ov{M}_+(1,t))$.  Note
that for $t>0$ sufficiently large, the surface $\ov{M}_+(1,t)$ consists of
two almost flat annular end representatives for $\ov{M}$ and
$\partial\ov{M}_+(1,t) = \partial \ov{M}_-(1,t)=\ov{M} \cap \ov{M}_-(1,t).$ 

Define $t_1$ to be the infinum of the values $t$ such that
$\overline{M}_+(1,t)$ is a nonnegative graph with bounded gradient
over its projection to $E(1,t)$ and $\partial \ov{M}_-(1,t) =\ov{M}
\cap \ov{M}_-(1,t)$.
 By the interior maximum principle and the Hopf maximum
principle applied along $\partial \ov{M}_-(1, t_1) $, we
observe that $\ov{M}_+(1,t_1)$
is a graph over its projection to $E(1,t_1),\; R(1,t_1)
(\ov{M})= \ov{M}$ and $\ov{M}$ is orthogonal to $E(1,t_1).$ This
observation is just the standard one that arises in the application of
the Alexandrov reflection argument, when one take into account the
maximum principle at infinity \cite{mr1} which guarantees that an end
of $\ov{M}_-(1,t_1)$ cannot be asymptotic to an end of $\ov{M}$ unless
$\ov{M}_-(1,t_1) \subset \ov{M}$.  

Note that $\partial \ov{M}_+(1,t_1)$ consists of a finite number of
simple closed curves in $E(1,t_1)$, since it has compact boundary.  
Furthermore, $\ov{M}_+(1,t_1)$ is a
planar domain with two ends, since it is a graph over a proper 
noncompact planar
domain in the annulus $E(1,t_1)$ with two ends.  Since the Euler 
characteristic $\chi(\ov{M})=2(1-\vk)-4=-2-2\vk,$ then
$\partial \ov{M}_+(1,t_1) = \ov{M}
\cap E(1,t_1)$ consists of $\vk+1$ simple closed curves.  Since these
curves are planar lines of curvature on $\ov{M}$ and $\ov{M}_+(1,t_1)$
is a graph over its projection to $E(1,t_1),$ then the simple closed curves in $\partial M_+(1,t)$ are
strictly convex curves bounding disks in $E(1,t_1)$.

Similarly, we have for some $t_2$ a related annulus $E(2,t_2)$
parallel to $A_{x_2}$, which is
an Alexandrov annulus of symmetry for $\ov{M}$.  After a fixed
translation of $\ov{M}$, we may assume that the circle $E(1,t_1) \cap
E(2,t_2)$ is $(0,0) \times \re/\mathbb{Z}\subset \re^2 \times \re/\mathbb{Z}$.
The proposition now follows. \qed

\section{The angle map $\theta \colon \mathcal{M}_{\vk} \rightarrow (0,
  \infty]$ is a local diffeomorphism.} \label{localuniqueness}

Our goal in this section is to prove the following result.

\begin{prop} \label{diffeo} For any component $C$ in $\mathcal{M}_{\vk}$,
  the image $\theta(C)$ is an open subset of $(0,
  \frac{\pi}{2}]$ and $\theta|_C \colon C \rightarrow \theta(C)$ is a
  diffeomorphism.
\end{prop}

The proof of Proposition \ref{diffeo} will depend on the following lemma.

\begin{lem} \label{locally rigid}
If $\ov{M} \in \mathcal{M}_{\vk}$ with $\theta (\ov{M})=
  \theta_0$, then $\ov{M}$ is locally rigid, i.e., there are no
  deformations of $\ov{M}$ though a family of minimal surfaces in
  $\mathcal{M}_{\vk}$ with the same angle $\theta_0$.  \end{lem}

Before we begin the proofs, we need to recall the \wei representation
and set some notation.
Recall that for a Riemann surface $\SR$ with a holomorphic
function $g$ and a holomorphic form $dh$ (not necessarily
exact, despite the notation), we may define a conformal branched minimal immersion 
via a map $F\colon\SR\to \Bbb{E}^3$ by

\begin{equation} \label{wei rep}
z\mapsto\ree\ngth\int^z_{\rho_0}\(\tfrac12\(g-\tfrac1g\)dh,
\tfrac i2\(g+\tfrac1g\)dh,dh\).
\end{equation}

For this surface, the function $g$ will be the Gauss map (postcomposed
with stereographic projection) and $dh$ will be the
complexified differential of the third coordinate in $\Bbb{E}^3$.
Conversely, given a conformal minimal immersion $F\colon\SR\to \Bbb{E}^3$ with 
Gauss map $g$ and complexified differential $dh$ of the third
coordinate, the surface may be represented by the expression
\eqref{wei rep}.  The induced metric on the minimal surface is given 
by 
\begin{equation} \label{ds-alpha}
ds_{F(\mathcal{R})}=\f{1}{{2}}\(|g|+\f1{|g|}\)|dh|;
\end{equation} 
thus a regular minimal surface will have zeroes of $dh$ of order $n$
balanced by poles or zeroes of $g$ of the same order.

The global problem for
producing minimal surfaces is that of well-definedness: analytic continuation
around a cycle must leave the map unchanged.
Thus we require
\begin{equation}
\begin{split}
\ree &\int_\g\tfrac12\(g-\tfrac1g\)dh
=\ree\int_\g\tfrac i2\(g+\tfrac1g\)dh\\
&=\ree\int_\g
dh=0
\end{split}
\end{equation}
for every cycle $\g\subset\SR$.

With this background in hand, we may begin the proof of
Lemma~\ref{locally rigid}. 

\proof Let $M$ be the lift of $\ov{M}$ to $\rth$. 
Let $g$ and $dh$ denote the \wei data of $M$ in the standard
notation. These combine to define the forms $gdh$ and $\frac1gdh$, which
we restrict to one of the fundamental domains of the surface with 
respect to the
$\BZ_2\oplus\BZ_2$ group of Alexandrov reflections guaranteed to exist
by Proposition \ref{symmetry}. 

\subsection{The shape of the developments of $|gdh|$ and $|\frac1gdh|$.}\label{zigzag}
The fundamental 
domain described above is planar, and
as the group elements act as isometries of the singular flat metrics
$|gdh|$ and $|\frac1gdh|$, the forms $gdh$ and $\frac1gdh$ develop
injectively to planar Euclidean domains, say $\ogup$ and $\ogdn$ 
bounded by a `zigzag' boundary, as we describe in just a moment.
Before we do that, however, we note that we already see that as the
fundamental domain is periodic with respect to a translation (in space), 
which is an isometry
of $|gdh|$ and $|\frac1gdh|$ as well, then the developed image of the domain
is invariant by a cyclic group of translational isometries (of the plane).

Much of our attention in this section will be focussed on the
boundaries of the developed domains $\ogup$ and $\ogdn$, as those
boundaries contain much of the geometry of the surface $M$.
Indeed, from the geometry of $M$, we recognize the basic shape of the zigzag
boundary of the developed image in $\mathbb{E}^2$. To draw 
this boundary, begin by drawing
an arc downwards for some distance at a slope of $-1$. Then make a left
turn and draw a segment upwards at a slope of +1.  Then draw downwards along a
segment of slope $-1$, etc. Continue drawing in this way $n=2\vk+2$
segments ($n$ is even and $\vk$ is the genus of $M$) and then repeat the
pattern indefinitely, both backwards and forwards.
This construction is meant to determine
the $gdh$ structure (the portion of the plane above the zigzag) for 
the surface $M$.

The $\frac{1}{g}dh$ structure is then determined by the requirement that
its periods should be conjugate to those of the $gdh$ structure, i.e.
$\int gdh=\overline{\int_{\g}\frac1gdh}$ for all cycles $\g\subset M$. So
we do this: on another plane, we 
draw another zigzag, so that the portion of the
plane above that zigzag will correspond to the $\frac{1}{g}dh$ structure.  The
conjugate period requirement is that we draw it as follows: we draw the
first arc at a slope of +1 and of exactly the same length as the first
segment on the first ($gdh$) zigzag. Then we draw the second arc at a slope of
$-1$ of exactly the same length as the second segment on the first zigzag.
We draw the third arc at slope +1 and of the same length as the third segment
of the first zigzag. We continue this construction for all $n$
segments and then repeat in both directions to guarantee
symmetry by an infinite group of translations (i.e. isometries
of the range $\Bbb{E}^2$ of development).

We assert that the developments $\ogup$ and $\ogdn$ of a fundamental
domain for the action of the
$\BZ_2\oplus\BZ_2$ group of Alexandrov reflections on $M$ have the
forms described above. In effect, we have to prove two 
statements to justify this:  first we need to show that the
boundary of those domains $\ogup$ and $\ogdn$ are piecewise straight,
and then we need to show that those straight edges meet at angles that
alternate between $\f{\pi}{2}$ and $\f{3}{2}\pi$.  For the first claim,
note that those Alexandrov reflections are
isometries of the flat singular metrics $|gdh|$ and $|\frac{1}{g}dh|$,
and that the boundary of the fundamental domain is fixed by the
isometry.  But as a smooth fixed set of an isometry is totally
geodesic, and the metrics $|gdh|$ and $|\frac{1}{g}dh|$ are
flat and smooth away from the the poles and zeroes of $g$, we see that
the smooth components of the boundary of the developed images
($\ogup$ and $\ogdn$)
of the fundamental domain of $M$ is bounded by straight lines.
Next observe that at the endpoints of those geodesic segments
(those endpoints corresponding to the $2\vk +2$ points where
$\ov{M}$ meets the intersection $A_{x_2} \cap A_{x_1}$ described in 
Proposition~\ref{symmetry}),
the forms $gdh$ and $\frac{1}{g}dh$ alternate between being regular
and having second order zeroes. (Also, one of those two forms 
has a double zero at such a point if and only if the other one of the
forms is regular there.)  Of course, the developed image of 
a holomorphic one-form with a zero of order $k$ has a cone point
with cone angle $2\pi(k+1)$; as these endpoints of the straight 
lines are fixed points of all four of the reflections, we see that
one-quarter of the cone-angles will be visible in one of the
fundamental domains.  Thus the boundaries will be composed of straight
lines, meeting at angles that alternate between $\f{\pi}{2}$ and
$\f{3}{2}\pi$, with angles at corresponding points of 
$\ogup$ and $\ogdn$ disagreeing, as claimed.

Because the $gdh$ and $\frac{1}{g}dh$ structures are defined on the
same Riemann surface, there is a conformal map between
those planar domains (above the corresponding zigzag). 

In all of this, we have ignored the flat structure for 
the form $dh$. This is because, following the same procedure
for the development of $dh$ as we did for the developments
of $\ogup$ and $\ogdn$, we see that the form $dh$ develops
on a fundamental domain to a domain with piecewise straight 
edges meeting at angles of $\pi$.  Thus the developed
image is a (geometric) halfplane with a periodic
collection of distinguished points on the boundary.  As any such
domain satisfies the (vertical) period condition that 
the distinguished points lie horizontally parallel to each other,
we see that any such domain will satisfy the relevant 
period condition, and there is no restriction on the
geometry of this domain.  This geometric fact corresponds to the
observation that the developed domain doubles
to a geometric sphere, so that the form $dh$ is exact.  From
both points of view, we cannot expect to glean much
information from the period condition for $M$ on $dh$.

\subsection{The angle between the ends.}

Having described the geometric structure of the developments
$\ogup$ and $\ogdn$ of the defining \wei
forms $gdh$ and $\frac{1}{g}dh$ of $M$, our next goal 
is to describe the moduli space $\SM(\theta_0)$ of
candidates for surfaces with asymptotic angle fixed at $\th=\th_0$.
To do this, we need to recognize the angle between those 
asymptotic planes in the
zigzags, as our moduli space will be defined in terms of zigzags.

\begin{prop} \label{prop6}
There is a function $m(\theta)$ which is strictly monotone
in $\theta$ so that if the Scherk ends of $M$ make an angle of 
$\theta$ with one another, then the zigzag boundaries of the 
domains $\ogup$ and $\ogdn$ are invariant by a group of translations
generated by $\<z\mapsto z+\ell_M(e)\exp(im(\th))\>$,
where $\ell_M(e)$ is the  length of the translation vector of $M$ and
$i = \sqrt{-1}$.
\end{prop}

We begin by considering the quotient of the surface by a single vertical
plane of symmetry. A single fundamental domain of the quotient has 
flux $\ora F$ across this
boundary equal to the flux across its ends, the latter given by
\begin{equation}\label{fluxeqn}
\ora F = (2\cos (\theta/2))h\vec n.
\end{equation} 
Here $\th$ is the angle between the planes, $h$ is the (normalized) height
of the fundamental domain, and $\vec n$ is the (appropriate) normal to the
reflective planes. Of course, the pair of ends of the domain is
homologous 
to the $\vk+1$
closed curves $\{\g_1,\dots,\g_{\vk+1}\}$ of intersection between the plane and
the surface. As these curves are orthogonal to the plane, the flux across
those circles is given by
$$
\begin{aligned}
\ora F  &= \sum^{\vk+1}_{i=1}\ell_M(\g_i)\vec n\\
&= \left(\sum^{\vk+1}_{i=1}\int_{\g_i}\frac{1}{2}(|g| +
 \frac1{|g|})|dh|\right)\vec{n}\\
&= \frac{1}{{2}}\left(\sum^{\vk+1}_{i=1}\left[\int_{\g_i}|gdh| +
  \int_{\g_i}|\frac1gdh|\right] \right)\vec{n},
\end{aligned}
$$
as the line element on the surface is given by
$ds_M=\frac{1}{2}(|g|+\frac1{|g|})|dh|$. Thus, we can rewrite the
length of the flux vector as
$$
F=|\ora F|= \frac{1}{2}\sum^{\vk+1}_{i=0}[\ell_{|gdh|}(\g_i) + \ell_{|gdh|}(\g_i)].
$$

We now consider the
quotient by action of reflection in 
the other vertical plane. As reflection about this plane
is an isometry for the metrics $|gdh|$ and $|\frac1gdh|$, we see that we
may rewrite the previous equation as 
$$
F = \sum^{\vk+1}_{i=0}[\ell_{\ogup}(\G_i) + \ell_{\ogdn}(\G_i)].
$$
Here, we multiply by two because we are only measuring lengths on a
single fundamental domain of the $\BZ_2\oplus\BZ_2$ action
(instead of in the pair of fundamental domains in the previous line), 
and we interpret
$\g_i$ as having trace $\G_i$ in that domain. Yet, by construction, the arcs
$\G_i$ are the arcs of the zigzags which all have the same slope, say $+1$.
Since, by  construction as well, we have
$$
\ell_{\ogup}(\G_i) = \ell_{\ogdn}(\G_i),
$$
we conclude
\begin{equation}\label{fluxlength}
F = 2\sum^{\vk+1}_{i=0}\ell_{\ogup}(\G_i).
\end{equation}

We consider next the total translational displacement of the zigzag,
i.e. the 
Euclidean distance between a point and its image under a generator of
the isometry group of $\ogup$ or $\ogdn$.  Elementary 
Euclidean geometry describes this quantity in terms of 
the the segments of the zigzag, as follows. Note that if we project 
a zigzag along one of its directions, we obtain a segment composed of 
isometric images of the the arcs of one slope, and if we project in
the orthogonal direction, we obtain a segment composed of isometric images
of arcs of the other slope: the translational displacement is the
length of the hypotenuse of the right triangle with these two segments 
as legs. Now the total translational displacement
is normalized by the
requirement that the translation $\ell_M(e)$
of the end is fixed; if $\G$ is an arc
around an end of a fundamental domain, then
$\ell_M(\G)=\frac{1}{2}(\ell_{|gdh|}(\G)+\ell_{|g^{-1}dh|}(\G))=\ell_{\ogup}(\G)$.
Thus, $\ell_M(e)= \ell_{\ogup}(\G) +o(1)$ (as $\G$ tends towards the 
end $e$)
and this fixes the translational displacement length . 

In summary, from equations
\eqref{fluxeqn} and \eqref{fluxlength}, we know that the total 
length, say $L_+(\th)$, 
in $\ogup$ of the segments $\{\G_0,\dots,\G_{\vk}\}$ (of slope $+1$)
is a monotone function of $\th \in (0, \pi)$, 
while the length of the total displacement vector
of the segments in a fundamental domain has length
fixed independently of the angle $\th$. Now the difference between
the endpoints of a fundamental domain of the zigzag is 
described as having length $\ell_M(e)$ and argument given as
$\arctan(\frac{L_+}{L_-}) - \frac{\pi}{4}$, 
where $L_-$ is the total length of all
of the segments of the zigzag of slope $-1$. Thus, as
$L_-^2 + L_+^2 = \ell_M(e)^2$, we see that the 
the slope of the orbit of a point in $\ogup$ under the
action of the cyclic group of translations is 
monotone in the asymptotic angle $\th$
between the ends. \qed

For the rest of the proof of Lemma \ref{locally rigid}, we will assume that the angle between the
ends is fixed.  The passage above shows that this forces the zigzag
boundary to have an orbit whose slope is $m(\theta_0)$, which is a well-defined
constant depending only on $\theta_0$.

\subsection{Deformations of Orthodisks.}\label{orthodisk deformation}

In general, a domain bounded
by a zigzag with orthogonal edges is known as an {\it orthodisk}.  An
orthodisk has geometry described by the positions of its vertices
$\{P_i\}$.  A pair of orthodisks with conjugate edge vectors
$\overrightarrow{P_jP_{j+1}}$ are called {\it conjugate} orthodisks. 
We have shown in section~\ref{zigzag} how a minimal surface of the
type we are considering in this paper gives rise to a conjugate
pair of orthodisks. Moreover, such a conjugate pair of orthodisks
is quite special, as the identity map on the minimal surface 
descends to a conformal map between the orthodisks which takes
vertices on one orthodisk, say $\ogup$, to corresponding
vertices on the other orthodisk, say $\ogdn$.  In this 
subsection, we will study a moduli space of pairs of 
conjugate orthodisks; these pairs will usually not be related
by a conformal map which preserves corresponding vertices.

To introduce this space, consider a surface $\ov{M} \in \SM_{\vk}$ with 
$\th(\ov{M}) = \theta_0$, as described in the statement of Proposition
\ref{prop6}. Then the domains $\ogup$ and $\ogdn$ for
$\ov{M}$ have zigzag boundaries with slopes $m(\th_0)$ as described 
in the previous subsection.  There is also then 
a $2\vk$-dimensional family $\{(\ogup, \ogdn)\}$ of pairs of domains
bounded by deformations of those zigzags invariant under the same
group of Euclidean planar isometries as for $\ov{M}$: of course,
most of the pairs in this space will not be related
by a conformal map which preserves corresponding vertices, as would
be the case for the pair, say $Z_0$, of domains for $\ov{M}$.

We then consider a family $\{M_t\}$ of minimal surfaces containing
$M=M_0$; these then induce, as above,  a family
$\{Z_t\}$ of such pairs of zigzags passing through $Z_0$, 
which would deform through domains with
zigzag boundaries.  Infinitesimally, then, we can compute the general form
of the (infinitesimal) Beltrami differential for the $gdh$ domain and for its
counterpart on the $\frac{1}{g}dh$ domain.  This pair represents a tangent
direction to the pair ($gdh, \frac{1}{g}dh$), construed to be a point in
the product of Teichmuller spaces of the  quotient  domains.

With all of this background, the crux of the argument is to compute
those Beltrami differentials. In particular, let us denote by, say
$\nu$, the Beltrami differential  prescribing the deformation on the $\frac{1}{g}dh$
structure, and by, say $\mu$, the Beltrami differential parametrizing the
deformation of the $gdh$ structure.  We then pull back $\nu$ to the $gdh$
structure via the assumed conformal map $F$ and obtain a Beltrami
differential
$F^*\nu$.  Then, if $\mu$ is the Beltrami differential for the $gdh$ structure,
we compute (!)   $F^*\nu = -\mu$.

The upshot of this result is that the two domains $\ogup$ and $\ogdn$
cannot remain conformally related for $t \not= 0$,
unless the lengths of all of the segments are preserved, 
ie., $gdh$ and $\frac{1}{g}dh$
do not change, 
which means that the family $\{M_t\}$ of minimal surfaces is infinitesimally
unmoving.

In particular, from  $F^*\nu = -\mu$,
and since $-\mu$ is not equivalent to $\mu$ unless both are equivalent to
zero, we conclude that neither structure has deformed.  But one can check
that this can only mean that no periods have changed, and so, the minimal
surfaces $\{M_t\}$  have 
only deformed by a congruence/homothety. That concludes
the argument, at least in outline form.

We need to formalize the previous discussion.  Let $\ogup$ and
$\ogdn$ denote the orthodisk structures for the forms $gdh$ and
$\frac1gdh$, respectively.  We are concerned with relating the
Euclidean geometry of the  orthodisks (which corresponds directly with the
periods of the
\wei data, as in the construction above) 
to the conformal data of the domains $\ogup$ and
$\ogdn$. From the discussion above, since a family of \mss\ $\{M_t\}$ will
always give rise to a corresponding family $\{(\ogup(t),\ogdn(t))\}$ of
orthodisks, it is clear that the allowable infinitesimal motions can be
parameterized in terms of the Euclidean geometry of 
$\ogup$ and $\ogdn$. These infinitesimal motions are given by infinitesimal changes in
lengths of finite sides  with
the changes being done simultaneously on $\ogup$ and $\ogdn$
to preserve conjugacy of the periods. The link to the conformal geometry is
that a motion which infinitesimally transforms 
$\ogup$, say, will produce an infinitesimal change in the 
conformal structure. Tensorially, this tangent vector to the
moduli space of conformal structures is represented by a
Beltrami differential. 

\subsubsection{Infinitesimal pushes.}

Here, we explicitly compute the effect of infinitesimal pushes
of  certain edges on the conformal geometry. 
This is done by explicitly displaying the infinitesimal
deformation and then computing the Beltrami differential.

In what follows, we rotate our picture by an angle of $\pi/4$ so that
all of our boundary edges are either horizontal or vertical.  This
simplifies our notation somewhat, even if it complicates the
meaning of 'conjugacy'.

\begin{figure}[H]
        \centerline{ 
               \includegraphics[width=5in]{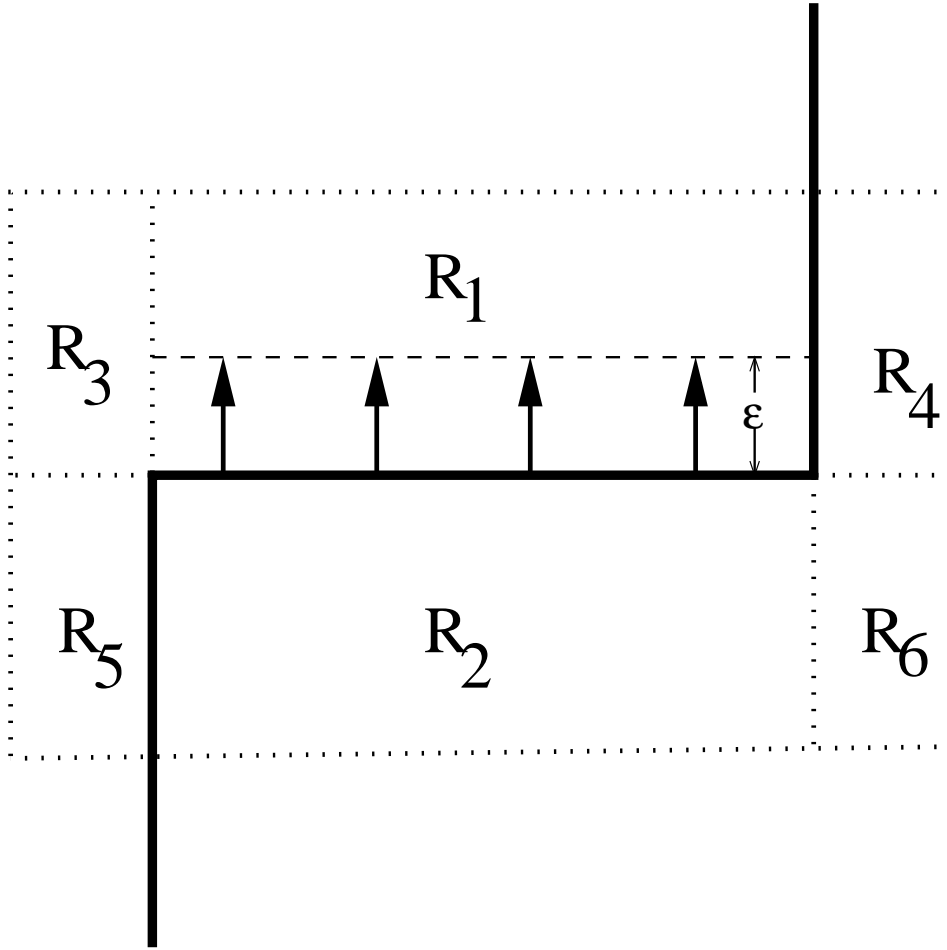}
      }
        \caption{Domain for Beltrami Differential Computation}
        \label{figure:deformation}
\end{figure}

We first consider the case of a horizontal finite side, say $E$; as in
the figure above, we see that the neighborhood of the
horizontal side of the orthodisk in the plane naturally
divides into six regions which we label $R_1$,...,$R_6$. Our
deformation $f_\ve = f_{\ve, b, \delta}$ differs from the identity
only in such a neighborhood,  and in each of the six regions,
the map is affine.
In fact, we have a two-parameter family of these deformations, all
of which have the same infinitesimal effect, with the parameters
$b$ and $\delta$ depending on the dimensions of the supporting neighborhood. 

\begin{equation}
f_\ve(x,y) = \begin{cases}
\(x,\ve +\frac{b-\ve}b y\),   &\{-a\le x\le a, 0\le y\le b\}=R_1\\
\(x,\ve +\frac{b+\ve}b y\),   &\{-a\le x\le a, -b\le y\le0\}=R_2\\
\(x,y+\frac{\ve+\frac{b-\ve}b y-y}\delta(x+\delta+a)\),   &\{-a-\delta\le x\le-a,
0\le y\le b\} = R_3\\
\(x,y-\frac{\ve+\frac{b-\ve}b y-y}\delta(x-\delta-a)\),   &\{a\le x\le a+\delta,
0\le y\le b\} = R_4\\
\(x,y+\frac{\ve+\frac{b+\ve}b y-y}\delta(x+\delta+a)\),   &\{-a-\delta\le x\le-a,
-b\le y\le0\} = R_5\\
\(x,y-\frac{\ve+\frac{b+\ve}b y-y}\delta(x-\delta-a)\),   &\{a\le x\le a+\delta,
-b\le y\le0\} = R_6\\
(x,y)   &\text{otherwise,} 
\end{cases} 
\end{equation}
where we have defined the regions $R_1,\dots,R_6$ within the
definition of $f_\ve$. Also, note that here, the orthodisk contains the arc
$\{(-a,y)\mid0\le y\le b\}\cup\{(x,0)\mid-a\le x\le a\}\cup
\{(a,y)\mid-b\le y\le0\}$.
Let $E$ denote the edge being pushed, defined above as $[-a,a]\times\{0\}$.

Let $\nu_\ve = \frac{\(f_\ve\)_{\bar z}}{\(f_\e\)_z}$
denote the Beltrami differential of $f_\ve$, and set
$\dot\nu=\frac d{d\ve}\bigm|_{\ve=0}\nu_\e$.  It is easy to compute that
$\dot\nu=[\frac d{d\ve}\bigm|_{\ve=0}\(f_\ve\)]_{\bar z}$ 
evaluates near $E$
to be
\begin{equation} \label{nu dot}
\dot\nu = \begin{cases}
\frac1{2b},   &z\in R_1\\
-\frac1{2b},   &z\in R_2\\
\frac1{2b}[x+\delta+a]/\delta+i\(1-y/b\)\frac1{2\delta}=\frac1{2b\delta}(\bar z+\delta+a+ib),
&z\in R_3\\
-\frac1{2b}[x-\delta-a]/\delta-i\(1-y/b\)\frac1{2\delta}=\frac1{2b\delta}(-\bar z+\delta+a-ib),
&z\in R_4\\
-\frac1{2b}[x+\delta+a]/\delta+i\(1+y/b\)\frac1{2\delta}=\frac1{2b\delta}(-\bar z-\delta-a+ib),
&z\in R_5\\
\frac1{2b}[x-\delta-a]/\delta-i\(1+y/b\)\frac1{2\delta}=\frac1{2b\delta}(\bar z-\delta-a-ib),
&z\in R_6\\
0  &z\notin\supp(f_\ve-\id).
\end{cases}
\end{equation}

\subsubsection{Deformation of a Conjugate Pair.}

With these definitions in place, we prove that a conformal pair
$(\ogup,\ogdn)$ of conjugate orthodisks 
(recall that this refers to a pair of zigzags with conjugate edge vectors)
admits only trivial deformations.
Let the space $\Delta_{\vk}$ denote the space of conjugate pairs of orthodisks
of the type described at the opening of the argument in 
subsection~\ref{orthodisk deformation}. The two-fold
branched cover of the double of one of these orthodisks is, after quotient
by a translation, a surface of genus $\vk$ 
(which branches over the sphere at a set of $2\vk +2$ points on an equator). 
Thus, we may regard $\Delta_{\vk}$ as
a subset $\Delta_{\vk}\subset \ST_{\vk}\x \ST_{\vk}$ of the 
product of a pair of \tec
spaces $\ST_{\vk}$; note
$\dim \Delta_{\vk}=\dim\ST_{\vk}$. 
Let $X_0$ denote a pair which is related by a conformal map
$\z:\ogup\to\ogdn$ (and which  therefore defines a periodic \ms\ with the
appropriate ends); note that such an example may be regarded 
as an element $X_0 \in \Delta_{\vk}\subset \ST_{\vk}\x \ST_{\vk}$.

We claim that $\Delta_{\vk}$ is transverse to the diagonal
$D_{\vk}=\{(\SR,\SR)\}\subset \ST_{\vk}\x \ST_{\vk}$,
where of course $\dim D_{\vk}=\dim \Delta_{\vk}$. To see this, 
note that we may regard a tangent direction as a pair
$(\dot{\nu}_{gdh},\dot{\nu}_{g^{-1}dh})$ of Beltrami differentials, each
representing a tangent direction to the points $[\ogup]\in
\ST_{\vk}$ and $[\ogdn]\in \ST_{\vk}$, respectively. Yet at
$X_0$, the points $[\ogup]$ and $[\ogdn]$ represent the
identical point in $\ST_{\vk}$, so we begin by computing how
the Beltrami differentials $\dot{\nu}_{gdh}$ and $\dot{\nu}_{g^{-1}dh}$
relate to one another. To this end, consider how an
infinitesimal push in the sense of the previous section on an edge $E$
defines Beltrami differentials $\dot{\nu}_{gdh}$ and
$\dot{\nu}_{g^{-1}dh}$. Of course, the conjugacy of $\ogup$ and
$\ogdn$ provides, via the formulas of that section, the basic
defining relation that if $\dot{\nu}_{gdh}$ has local expansion
$\dot{\nu}_{gdh}(z)=\frac1{2b}\frac{d\bar z}{dz}$ near an interior point
of an edge $E$, then also
$\dot{\nu}_{g^{-1}dh}(\z)=\frac1{2b}\frac{d\bar\z}{d\z}$
at the corresponding edge. However, since
$X_0$ is a pair of conformally related orthodisks, we may also assume, in
this particular case, the existence of a conformal map $\z:\ogup\lra\ogdn$
which preserves the vertices. Such a map takes vertical sides
to horizontal sides by construction and this has the local
expansion $\z=\pm i|c|z+0(|z|^2)$ near an interior point of an
edge. We therefore compute the pullback of $\dot{\nu}_{g^{-1}dh}$
to $\ogup$ as
$$
\begin{aligned}
\dot{\nu}_{g^{-1}dh}(\z)\frac{d\bar\z}{d\z} &= \dot{\nu}_{g^{-1}dh}(\z)
\frac{\bar\z'}{\z'}\frac{d\bar z}{dz}\\
&= (-\dot{\nu}_{g^{-1}dh}(\z) + h.o.t.)d\bar z/dz\\
&= -\frac1{2b}\frac{d\bar z}{dz}
\end{aligned}
$$
along the edge $E$. 
We conclude that, locally in the interior of the edge $E$, we have
$\z^*\dot{\nu}_{g^{-1}dh}=-\dot{\nu}_{gdh} +o_b(1)$, where $o_b(1)$ may
be taken arbitrarily small. The above computation nearly proves that
$[\dot{\nu}_{gdh}]=-[\dot{\nu}_{g^{-1}dh}]$ as elements of the tangent
space $T_{\ogup} \ST_{\vk}=T_{\ogdn} \ST_{\vk}$: what is left to prove
for that conclusion is that the contributions to 
$[\dot{\nu}_{gdh}]$ or $[\dot{\nu}_{g^{-1}dh}]$ from the regions
$R_3$ and $R_5$ -- and their counterparts in the other 
fundamental domains for the 
$\BZ_2\oplus\BZ_2$ group of Alexandrov reflections -- is
negligible.  To see this, begin by noticing that our surface  
$\ov{M}$ is hyperelliptic, branching over the points where
the Gauss map $g$ is singular, ie. over the images of the
endpoints of the edges of the zigzags. Recall next that the pairing of the 
tangent and cotangent spaces to $\ST_{\vk}$ is given by integrating
the product of Beltrami differentials and holomorphic quadratic
differentials on 
$\ov{M}$. Now, as $\ov{M}$ is hyperelliptic, the Taylor 
expansion of a holomorphic quadratic differential
$\Phi$ on $\ov{M}$ around (the lift of) a corner of an edge is even,
i.e. $\Phi = (c_0 +c_2w^2 + h.o.t)dw^2$.  This descends to 
a neighborhood of a corner of the edge via the
local map $w=z^{\frac{1}{3}}$ as 
$\Phi = (\frac{1}{9}c_0z^{-\frac{4}{3}} + \frac{1}{9}c_2z^{-\frac{2}{3}}
+ O_z(1))dz^2$.  

In terms of this expression, the terms of
order $O(z^{-\frac{2}{3}})$ and higher are easily seen to 
integrate against $[\dot{\nu}_{gdh}]$ or $[\dot{\nu}_{g^{-1}dh}]$ to  
negligible amounts in a small neighborhood of the endpoints of
an edge, but the term of order $O(z^{-\frac{4}{3}})$ is more worrying.
To understand that it also has but a negligible effect, we
need to consider its effect in an entire neighborhood in
$\ov{M}$ of an endpoint of an edge of the zigaag.
In particular, we consider the development of the 
other three fundamental domains whose closures contain that
point: these developments are obtained by reflection 
across the edges incident to that point.  After a
reflection in each of these edges, we obtain the development 
of the domain in the image of $w \mapsto -w$ of our original
domain:  it has the form in Figure~\ref{figure:deformation} consisting
of the regions $R_2, R_4$ and $R_6$.  In that region,
because the deformation of the relevant edge will be
downwards instead of upwards, the Beltrami
differentials $\dot{\nu}_{g^{-1}dh}$ and $\dot{\nu}_{gdh}$
should be regarded as expressed by the 
formulas for $-\dot{\nu}$ in \eqref{nu dot}.
In terms of these explicit formulas and using that
the map between these developments may be written
as $\Psi: z \mapsto -z$ (abusing notation by setting $a=0$ in both  
formulas), we note that 
\begin{equation}
\Psi^*(-\dot{\nu}\big|_{R_6}) + \dot{\nu}\big|_{R_3} = 
\Psi^*(-\dot{\nu}\big|_{R_4}) + \dot{\nu}\big|_{R_5} = 0.
\end{equation}
On the other hand, since $\Phi$ is even, we have that 
$\Psi^*\Phi =\Phi$, and so we conclude that the 
pairing $\int \Phi\dot{\nu}$
between $\Phi$ and either $\dot{\nu}=\dot{\nu}_{g^{-1}dh}$ or 
$\dot{\nu}=\dot{\nu}_{gdh}$
over the images of the regions $R_3, R_4, R_5$ and $R_6$
(in all of the fundamental domains) vanishes: the effect
of  $\dot{\nu}_{g^{-1}dh}$ and $\dot{\nu}_{gdh}$ as tangent vectors
in a pairing with covectors is given by integration over
(the images of) the regions $R_1$ (and $R_2$).

We conclude that for $\dot{\nu}_{g^{-1}dh}$ and $\dot{\nu}_{gdh}$
defining an infinitesimal push, we have
$[\z^*\dot{\nu}_{g^{-1}dh}]=-[\dot{\nu}_{gdh}]$ as elements of the
tangent space to $\ST_{\vk}$. Further, 
since any deformation of $X_0$ is given by a linear
combination of such infinitesimal pushes, we conclude that
$[\dot{\nu}_{gdh}]=-[\dot{\nu}_{g^{-1}dh}]$ as elements of the tangent
space $T_{\ogup} \ST_{\vk}=T_{\ogdn} \ST_{\vk}$. Thus, any 
conjugacy-preserving deformation of $X_0$ destroys the conformal
equivalence of $\ogup$ and $\ogdn$ to the order of the
deformation, a statement which implies $\Delta_{\vk}$ is transverse to the
diagonal. This concludes the proof of the claim.

To finish the proof of Proposition \ref{diffeo}, observe any deformation $M_t$ 
of minimal surfaces through $M_0$ must
preserve the conformality between $\Omega_{gdh}(t)$ and
$\Omega_{g^{-1}dh}(t)$.  Thus, by the computation above, we conclude that the
tangent vector for
$M_t$ is trivial as a tangent vector to \tec space, and moreover, the
forms $gdh$ and
$\frac1gdh$ are only trivially deformed. Since these forms 
suffice as Weierstrass
data, we conclude that $M_t$ deforms only by an infinitesimal congruence,
as desired. 

Recall the moduli spaces $\mathcal{M}(\theta)$ of pairs of surfaces
$(\Omega_{gdh}, \Omega_{g^{-1}dh})$ whose zigzag boundaries have
orbits of points which accend at slope $m(\theta)$ under the
translation group.  The paragraph above shows that
$\mathcal{M}(\theta_0)$ meets the 
diagonal $D_{\vk}$ in $\ST_{\vk} \times \ST_{\vk}$
transversely.  As the $\dim (\mathcal{M}(\theta)) + \dim (\ST_{\vk}) = \dim
(\ST_{\vk} \times \ST_{\vk})$, the 
implicit function theorem implies that there is
a curve $X_\theta \in \mathcal{M} (\theta) \cap D_{\vk}$ for $\theta$ near
$\theta_0$.  By our construction of zigzags, such an element
$X_\theta$ represents a periodic minimal surface whose ends make an
angle $\theta$ with each other.  This concludes the proof of
Proposition \ref{diffeo}.  \qed

\section{Properness of the angle map.}

In this section, we prove the following properness result.

\begin{prop} \label{proper} The angle map 
$\theta \colon \mathcal{M}_{\vk} \rightarrow
  (0, \frac{\pi}{2}]$ is proper.
\end{prop}

\proof Since $\mathcal{M}_{\vk}$ consists of curves and $\theta \colon \mathcal{M}_{\vk} \rightarrow (0,
\frac{\pi}{2}]$ is a local diffeomorphism by Proposition \ref{diffeo},
it is sufficient to prove that if $\{\ov{M}_n\}_{n \in \mathbb{N}}$ is
  a sequence of examples in $\mathcal{M}_{\vk}$ with $\lim_{n \rightarrow
    \infty} \theta (\ov{M}_n)= \theta_0>0$, then a subsequence of the
  surfaces converges on compact subsets of $\rth/\mathbb{Z}$ to a
  minimal surface $\ov{M}_\infty \in \mathcal{M}_{\vk}$ with
  $\theta(\ov{M}_\infty)= \theta_0$.  Note that we choose the surfaces
  $\{\ov{M}_n\}_{n \in \mathbb{N}}$ so that each lies in
  $\rth/\mathbb{Z}$ and is invariant under reflection in the vertical
  annuli $A_{x_1}, A_{x_2}$ in $\rth/\mathbb{Z}$ given in the
  statement of Proposition \ref{symmetry}.

As the four ends of $\ov{M}_n$ are asymptotically flat, the lifted
surface $M_n$ in $\rth$ has area growth constant $A(M_n)=2\pi$.  Thus,
by the monotonicity formula for area of a minimal surface in $\rth$,
every ${M}_n$ has at most $2 \pi r^2$  of
area in balls of radius $r$.  Hence, after choosing $r<\frac{1}{2}$ 
so that the balls in $\rth$ embed in the quotient space, we see that the 
surfaces $\ov{M}_n$ have locally bounded area in $\rth/\mathbb{Z}$.  We now
check that there are also uniform estimates for the curvature of the
surfaces in the sequence.

Arguing by contradiction and after extracting a subsequence, suppose that
there exist points $p_n \in \ov{M}_n$ with maximal absolute curvature
$\lambda_n^2 \geq n$; note that maximal curvature points $p_n$ exist since
the asymptotic curvature of the $\ov{M}_n$ is zero.  Consider the
homothetically expanded surfaces $\widetilde{M}_n=
\lambda_n[\ov{M}_n \cap B(p_n, \frac{1}{2})]$ in the homothetically expanded ball $\lambda_n B(p_n, \frac{1}{2}).$
The surfaces $\widetilde{M}_n$ are submanifolds
 in the flat
three-manifolds $\lambda_n B( p_n, \frac{1}{2})$, which are isometric
to balls $B(\vec{0},r_n)$ of radius $r_n=\frac{\lambda_n}{2} \geq \frac{\sqrt{n}}{2}$ in $\rth$
centered at the origin $\vec{0}$; these balls are converging naturally
to all of
$\rth$.  Consider the surfaces $\widetilde{M}_n$ to lie in $B(\vec{0},
r_n) \subset\rth$.
Since the $\widetilde{M}_n$ have maximal absolute curvature 1 at the origin
and in balls of radius $r \leq r_n$ have area at most $2\pi r^2$,
standard results (see, for example, \cite{mpe1}) imply that a
subsequence of these surfaces converges on compact subsets of $\rth$
to a properly embedded minimal
surface $\widetilde{M}$ in $\rth$ with absolute curvature at most $1$
and with absolute curvature $1$ at the origin.  
The surface $\widetilde{M}$ is connected by the strong halfspace
theorem \cite{hm10}. Furthermore, since each surface
$\widetilde{M}_n$ has the same total absolute curvature as $\ov{M}_n
\cap B(p_n, \frac{1}{2}),$ then each $\widetilde{M}_n$
has total absolute curvature less than the total absolute
curvature of $\ov{M}_n$, which by the Gauss-Bonnet formula is $-2 \pi
\chi (\ov{M}_n)$ which is the finite number $4 \pi(\vk+1)$.  In
particular, $\widetilde{M}$ has finite total curvature.  Since
$\widetilde{M}$ is
embedded and not flat, it has at least two ends which are asymptotic
to either planes or ends of catenoids.  Since $\widetilde{M}$ has area
growth constant at most $2\pi$ and it is not a plane, then it has exactly two ends.  By
Schoen's theorem \cite{sc1}, we see that $\widetilde{M}$ is a catenoid with waist
circle passing through the origin.

We claim that each of the Alexandrov annuli of symmetry of the
$\ov{M}_n$ intersect $B(p_n, r_n)$ for $n$ large and limit to
planes of symmetry for $\widetilde{M}$.  Otherwise, 
the surface $\widetilde{M}$ would be
the limit of domains in $\widetilde{M}_n\subset\lambda_n \ov{M}_n$, which are
graphical over their projections to one of its Alexandrov annuli of
symmetry.  It would then follow that the Gaussian image of $\widetilde{M}$
would lie in a hemisphere of $S^2$, which is false for a catenoid.
Hence, the Alexandrov annuli of symmetry of the $\widetilde{M}_n$ limit in
a natural way to Alexandrov planes of symmetry of $\widetilde{M}$, which,
after a translation of $\widetilde{M}$, we can consider to be the $(x_1,
x_3)$- and $(x_2, x_3)$-planes.  Furthermore, since each component
of the fixed point set of an Alexandrov annulus of symmetry is a
simple closed convex curve, which is invariant under reflection across
the other such Alexandrov annulus and hence has two fixed points, then it
is easy to show that, since $\widetilde{M}$ is connected, the fixed
point set of one of the planes of Alexandrov symmetry of $\widetilde{M}$
must intersect the fixed point set of the other plane of symmetry of
$\widetilde{M}$.  It follows that $\widetilde{M}$ is a catenoid
with axis being either the $x_1$-axis or the $x_2$-axis (rather than
with axis being the $x_3$-axis).

From the discussion in the previous paragraph, we conclude that a subsequence
of the locally finite integral minimal varifolds $\ov{M}_n$ in $\rth/\mathbb{Z}$ converges to
a limit minimal varifold $\ov{M}_\infty$ with mass density $2$ at some point
of the vertical circle $\alpha=(0,0) \times \re/\mathbb{Z}$ in $\rth/\mathbb{Z}$. 
It follows from the monotonicity formula for area that $\ov{M}_\infty$ is the union of
two flat vertical annuli, not necessarily distinct and both containing
$\alpha$.

We claim that the convergence of the $\ov{M}_n$ to $\ov{M}_\infty$ is smooth
away from $\alpha$.  If not, then 
there exists a point $p \in \rth/\mathbb{Z} - \alpha,$ such that, 
after extracting a subsequence, the absolute
curvature of the $\ov{M}_n$ in the $\ve = \frac{1}{2} d(p, \alpha)$
ball $B(p, \ve)$ centered at $p$ is at least $n$.  Let $q_n$ be a
point of $\widehat{M}_n =\ov{M}_n \cap B(p, \ve)$, where the function
$d(\cdot , \partial \widehat{M}_n) |K|(\cdot)$ has its maximum value;
here, $|K|(\cdot)$ is the absolute curvature function on $\widehat{M}_n.$  Let $\lambda_n =
\sqrt{|K(q_n)|}$ and note that the surfaces
  $\Sigma_n=\lambda_n(\widehat{M}_n \cap B(q_n, r_n)),$ where $r_n=
  \frac{1}{2} d(q_n, \partial B(p, \ve)),$ have bounded curvature in
  the balls $\lambda_n B(q_n,r_n)$ of radius $r_n$ centered at $q_n$.  These balls converge
  to $\rth$ with $q_n$ considered to be at the origin.  The
  surfaces $\Sigma_n$ have absolute curvature bounded by $4$
  and have local area estimates. Our previous arguments now imply
  that a subsequence of the $\Sigma_n$ converges to a catenoid in $\rth$ and this
  catenoid contains a point which is a limit of points coming
  originally from $\alpha=A_{x_1} \cap A_{x_2}$.  This is a contradiction, since $\alpha$
  is disjoint from $B(p, \ve)$.  This contradiction proves that the surfaces
  $\ov{M}_n$ converge smoothly with multiplicity two to $\ov{M}_\infty$, away
  from $\alpha$.

We claim that $\ov{M}_\infty$ is either $A_{x_1}$ or $A_{x_2}$ with
multiplicity two.  If not, then 
since  $\ov{M}_\infty$ is the union of two flat vertical annuli,
the circle $\alpha \subset
\ov{M}_\infty$ is contained in the intersection set of two distinct
vertical flat annuli, and so, every point $p \in \alpha$ is a point 
in the singular set
of convergence to $\ov{M}_\infty$.  
Now, the blow up argument in the
  previous paragraph shows that, for any $\ve >0$ the ball $B(p, \ve)$
  contains for $n$ large,
  an approximately scaled down catenoid in $\ov{M}_n$, and so, the 
total absolute
  curvature of $\ov{M}_n \cap B(p, \ve)$ is at least $3 \pi$ for $n$
large.  
Since
  $\ve$ is arbitrary and the total absolute curvature of $\ov{M}_n$ is
  $4 \pi (\vk+1)$ which is finite, we obtain a contradiction. Hence,
  $\ov{M}_\infty$ is $A_{x_1}$ or $A_{x_2}$ with multiplicity two as a
  limit varifold.  (Although we do not use it here, we observe that
this argument also
  shows, after choosing a subsequence, that there are at most $\vk+1$
  distinct singular points of convergence of the $\ov{M}_n$ to
  $A_{x_1}$ or $A_{x_2}$, which by our earlier arguments must lie on $\alpha$).

Let $N_\ve(\alpha)$ be any fixed $\ve >0$ neighborhood around
$\alpha$.  For $n$ large, our analysis of the limits $\lambda_n
\ov{M}_n$ shows that $M_n - N_\ve(\alpha)$ consists of four
annular Scherk ends of $M_n$.  Now, we have also shown that
the convergence of $\ov{M}_n- N_\ve(\alpha)$ 
to $\ov{M}_{\infty}-N_\ve(\alpha)$ is
smooth along $\partial (\ov{M}_n-N_\ve (\alpha))$, so the flux vectors
(integrals of the conormals along each component of $\partial
(\ov{M}_n- N_\ve (\alpha)))$ are converging to either $(\pm 1, 0,0)$
or $(0,\pm1, 0)$, since those limits are the 
flux vectors of $\ov{M}_{\infty}$.  But these
flux vectors are also the flux vectors of the ends of
$\ov{M}_n$ by the divergence theorem, and the flux vectors of the ends
of $\ov{M}_n$
are bounded away from $(\pm 1, 0,0)$, since the limit angle $\theta_0
>0$.  This contradiction implies that the sequence $\{\ov{M}_n\}_{n
  \in \mathbb{N}}$ with $\theta(\ov{M}_n) \rightarrow \theta_0 >0$
have uniformly bounded curvature.

We now prove that our original sequence $\{\ov{M}_n\}_{n \in
  \mathbb{N}}$ with $\theta (\ov{M}_n) \rightarrow \theta_0 >0$
  converges to an example $\ov{M}_\infty \in \mathcal{M}_{\vk}$.  Since
  the sequence of surfaces has uniformly bounded curvature and local
  area estimates, a subsequence converges on compact subsets of $\rth
  /\mathbb{Z}$ to a properly embedded minimal surface
  $\ov{M}_\infty$.  Recall that the sum of the lengths of the convex curves
  in $M_n \cap A_{x_i}$ corresponds to the flux of $\nabla x_{(i+1) \mod
  2}$ of
  $M_n$, which is less than or equal to 2 and is determined by
  $\theta (M_n)$. 
Since every convex curve in $\ov{M}_n
\cap (A_{x_1} \cup A_{x_2})$ intersects the vertical circle $\alpha$,
these $\vk+1$ convex curves $\ov{M}_n \cap A_{x_i}$ converge smoothly to
$\vk+1$ convex curves in $\ov{M}_\infty \cap A_{x_i}$, each of length
less than 2 for $i=1,2$. Moreover, the lengths
of these curves are also bounded away from zero, since they are
planar curves and principal on $M_n$: any pinching of them would
then blow up the curvature somewhere along them. Yet these 
lengths correspond to the lengths of the segments in the zigzags bounding
the domains $\ogup$ and $\ogdn$ for $M_n$, so we see that
these domains $\ogup$ and $\ogdn$ for $M_n$ also converge 
smoothly and without degeneration to the orthodisks of
$\ov{M}_\infty$.

These orthodisks, together with the implied vertex-preserving
conformal map between them, of course determine the \wei
data for a minimal surface whose geometry is given by the 
Euclidean geometry of the orthodisks.  Here, since the limiting
orthodisk has a fundamental domain bounded by $2\vk +2$ nondegenerate
segments of alternating slope, we see that the surface
$\ov{M}_\infty$ is a nondegenerate minimal surface of genus $\vk$.
Further, as the flux is determined (see Proposition~\ref{prop6})
by the slope of the orbit of a vertex, and the orthodisks
are converging smoothly, we see that $\theta(\ov{M}_\infty)=\theta_0$.

This completes the proof of the proposition.  \qed

\begin{rem} \label{rem8} We note that the flux argument given in the curvature
  estimate part of the
  proof of Proposition \ref{proper} implies that if, for some sequence in
  $\{\ov{M}_n\}_{n \in \mathbb{N}}$ in $\mathcal{M}_{\vk}$ we
have $\theta(\ov{M}_n)
  \rightarrow 0$, then the locally finite limit integral minimal  varifold
  $\ov{M}_\infty$ is the annulus $A_{x_1}$ with multiplicity two. 
\end{rem}

\section{Small angle examples are Scherk examples.}

In this section, we prove the following result.

\begin{prop}  \label{small} For every $\vk>0$, there exists an
  $\ve >0$ such that if $\ov{M} \in \mathcal{M}_{\vk}$ and $\theta
  (\ov{M}) < \ve$, then $\ov{M}$ is a Scherk example.
\end{prop}

\proof Suppose $\ov{M}_n \in \mathcal{M}_{\vk}$
is a sequence of examples, where $\theta (\ov{M}_n) < \frac{1}{n}$.
After extracting a subsequence, the $\ov{M}_n$ converges to an
integral varifold $\ov{M}_\infty$.  From the proof of Proposition
\ref{proper} and Remark \ref{rem8}, it is easy to see that $\ov{M}_\infty$ is one of the annuli
$A_{x_1}$ or $A_{x_2}$ 
(with multiplicity two) of symmetry and that the limiting flux vectors
to the ends of the $\ov{M}_n$ converge to vectors in $\ov{M}_\infty$.
Hence,  $\ov{M}_\infty$ corresponds to $A_{x_1}$.  Modifications
of the arguments used in the proof of Proposition \ref{proper} also
show that, for $\alpha=(0,0) \times \re/\mathbb{Z}$ and $n$ large, there exist $\vk+1$ points $P_n=\{p_1(n),
p_2(n), \hdots, p_{\vk+1}(n)\} \subset \ov{M}_n \cap \alpha$ with
normal vector $(0,0,1)$ together with small
positive numbers $\ve_1(n), \ve_2(n), \hdots, \ve_{\vk+1}(n)$, such that
for each $i$, the intersection
$B(p_i(n), \ve_i(n)) \cap \ov{M}_n$ is a compact annulus
which is $C^2$-close to a standard catenoid with axis along the $x_2$-axis,
which has been scaled by inverse of square root of the absolute
curvature at $p_i(n)$.  Furthermore, after replacing by a subsequence,
the surface $\widetilde{M}_n= \ov{M}_n \cap [\rth/\mathbb{Z} -
\bigcup^{\vk+1}_{i=1}B(p_i(n), \ve_i(n))]$ consists of two components which
are graphs of gradient less than $\frac{1}{n}$ over their 
projections to the annulus
$\ov{M}_{\infty}$.

A subsequence of the (paired) graphs $\widetilde{M}_n$ converges smoothly to
$\ov{M}_\infty$ punctured in at most $\vk+1$ points with graphical gradients
converging to zero as $n \rightarrow \infty$.  Thus, the degenerating
conformal structures of $\ov{M}_n$ as $n \rightarrow \infty$ can be seen to
converge to that of two copies of the annulus with nodes forming at the
(limits of the) $\vk+1$ points $P_n$ along $\alpha \subset \ov{M}_\infty$.

Now let $\wh A_n$ denote 
the $\Bbb{Z}$-cover of one of the fundamental annuli of
$M_n$ and let $\wh A$ denote the limit of $A_n$. Thus, the domains
$\wh A_n$ and $A_n$ are conformal half-planes, with the points $P_n$
lifting to a periodic sequence $\wh P_n$ of boundary points of $\wh
A_n$; that sequence $\wh P_n$ converges 
to a periodic sequence $\wh N_n$ of lifts  of the
nodes $N_n$ on the boundary of $\wh A$. We shall also have need of the
periodic sequence $\wh Q_n\subset\p\wh A_n$ which are lifts of the
points $Q_n\in\ov M_n\cap\a$ whose normal vectors are $(0,0,-1)$.
Naturally, the points $\wh P_n$ and $\wh Q_n$ alternate in position
along $\p\wh A_n$ and $\wh Q_n\to \wh N_n$ along with $\wh P_n$. We take the
upper half-plane $\BH$ as a model for $\wh A_n$ and $\wh A$, and we let
the images of $\wh P_n$ be given by $\{a_{k,\ell}=k+a_\ell\mid\ell=0,\dots,\vk$
where $a_0=0$ and $0<a_\ell<1$ for $\ell>1\}$ and the images of $\wh
Q_n$ be given by $\{b_{k,\ell}=k+b_\ell\mid\ell=0,\dots,\vk$, where
$0<b_\ell\le 1\}$. Naturally we take
$a_{k,j}<b_{k,j}<a_{k,j+1}<b_{k,j+1}$ for every $k$. These points
$a_{k,j}$ and $b_{k,j}$ depend on $n$, but we will suppress the natural
dependence on $n$ until it is relevant and important.

As in section~\ref{localuniqueness}, we consider the forms $gdh$ and
$\f1gdh$ on $\wh A_n$ and $A_n$. These evidently develop to domains
bounded by a periodic boundary, as described in
\S~\ref{zigzag}. We observe that we may parametrize the
domains $\ogup$ and $\ogdn$ via the Schwarz-Christoffel maps
(suppressing the dependence on $n$)
\begin{equation*}
F_{\ogup}(\z) = e^{-i\pi/4}\int^\z\prod^{\vk}_{j=0}
\[\f{z-a_{0,j}}{z-b_{0,j}}\]^{1/2}\prod^\infty_{k=1}\prod^{\vk}_{j=0}
\[\f{(z-a_{k,j})(z-a_{-k,j})}{(z-b_{k,j})(z-b_{-k,j})}\]^{1/2}dz
\end{equation*}
and
\begin{equation}\label{sc}
F_{\ogdn}(\z) = e^{-i\pi/4}\int^\z\prod^{\vk}_{z=0}
\[\f{z-b_{0,j}}{z-a_{0,j}}\]^{1/2}\prod^\infty_{k=1}\prod^{\vk}_{j=0}
\[\f{(z-b_{k,j})(z-b_{-k,j})}{(z-a_{k,j})(z-a_{-k,j})}\]^{1/2}dz.
\end{equation}
To see this, first observe that as the terms
$\(\f{(z-a_{k,j})(z-a_{-k,j})}{(z-b_{k,j})(z-b_{-k,j})}\)^{1/2}$ are
asymptotically $1+O(\f1{k^2})$ for $k$ large, the infinite product
converges absolutely, and uniformly on compacta in $\ov\BH$. Moreover,
the images of the boundary are evidently periodic zigzags: they are
zigzags by the basic Schwarz-Christoffel theory, and they are periodic
as the periodicity of the coefficients ($a_{k+1,j}=a_{k,j}+1$,
$b_{k+1,j}=b_{k,j}+1$) forces the periodicity of the developing maps
$F_{\ogup}$ and $F_{\ogdn}$.

Now, the crucial part of the analysis is the determination of the
coefficients $a_{k,j}$ and $b_{k,j}$: we know that as $n\to\infty$, we
have $|a_{k,j}-b_{k,j}|\to0$ so that $a_{k,j}$, $b_{k,j}\to c_{k,j}$,
and we need to determine both $|a_{k,j}-b_{k,j}|$ and $c_{k,j}$. This
sort of analysis has been carried out by Traizet in a number of
slightly different settings (e.g. \cite{tra1}, \cite{tra4}).
Unfortunately for the brevity of this argument, 
while we can follow his general outline, he does not
seem to have treated this precise case; fortunately, as all of the
relevant information about the surfaces $M_n$
is in the development \eqref{sc} (recall that the
$dh$ development offers no substantive restrictions), we can give a
full yet more elementary treatment directly from the equations
\eqref{sc}.

The crucial condition is that $\int_\g gdh=\int_\g\ov{\f1gdh}$ for
every cycle $\g\subset M$. On the annuli $A_n$, this implies that
$$
F_{gdh}(\a) - F_{gdh}(\b) = \ov{F_{g^{-1}dh}(\a) - F_{g^{-1}dh}(\b)}
$$
for $\a$, $\b\in P_n\cup Q_n$: here the point is that any cycle on
$M_n$ is homologous to a linear combination of arcs on the boundary
connecting the vertical points of the Gauss map.

Let us normalize the setting. We focus on four consecutive points
$\a-c\ve$, $\a$, $\b$, $\b+d\ve$; naturally, each of the points depends
on the parameter $n$, and as $n\to\infty$, the points $\a-c\ve$ and $\a$
converge to a node as do $\b$ and $\b+d\ve$. In the natural notation, we
compute
\begin{equation*}
\begin{aligned}
F_{gdh}(\b) - F_{gdh}(\a) &= c_0 + \f{e^{-i\pi/4}}2(d-c)\ve\log\ve +
\text{ higher order terms and}\\
F_{g^{-1}dh}(\b) - F_{g^{-1}dh}(\a) &= \ov{c_0} -
\f{e^{i\pi/4}}2(d-c)\ve\log\e + \text{ higher order terms}.
\end{aligned}
\end{equation*}
Thus, in order that
$F_{gdh}(\b)-F_{gdh}(\a)=\ov{F_{g^{-1}dh}(\b)-F_{g^{-1}dh}(\a)}$, we
must have (from the singular term) that $d=c$. As this computation
holds for the interval between any pair of points coalescing to a node,
we find that $|a_{k,j}-b_{k,j}|/|a_{k',j}-b_{k',j}|\to1$ as $n\to\infty$
for any choice of $k$, $k'$, $j$ and $j'$.

\begin{rem} This last statement reflects that the  sizes of the
curves (as curves in space) being pinched are (asymptotically) identical.

Before turning our attention to the ``small'' intervals between
$a_{k,j}$ and $b_{k,j}$, we readjust our notation, setting
$b_{k,j}=a_{k,j}+\ve+\eta_{k,j}$ where $\eta_{k,j}=o(\ve(n))=o(\ve)$ as
$n\to\infty$. Also, being mindful of convergence issues that will
eventually arise, we explicitly consider approximations
\begin{equation}
\begin{aligned}\label{sc-m}
F^M_{\ogup}(\z)  &= e^{i\pi/4}\int^\z\prod^{\vk}_{j=0}
\[\f{z-a_{0,j}}{z-b_{0,j}}\]^{1/2}
\prod^M_{k=1}\prod^{\vk}_{j=0}
\[\f{(z-a_{k,j})(z-a_{-k,j})}{(z-b_{k,j})(z-b_{-k,j})}\]^{1/2}dz\\
F^M_{\ogdn}(\z)  &= -e^{i\pi/4}\int^\z
\prod\[\f{z-b_{0,j}}{z-a_{0,j}}\]^{1/2}
\prod^M_{k=1}\prod^{\vk}_{j=0} \[\f{(z-b_{k,j})(z-b_{-k,j})}{(z-a_{k,j})(z-a_{-k,j})}\]^{1/2}dz
\end{aligned}
\end{equation}
to the infinite vertex Schwarz-Christoffel map in \eqref{sc}. As the
convergence of the ``partial product'' map in \eqref{sc-m} is uniform in
$M$ with estimates independent of $n$ (because $P_n\cup Q_n$ converge
to the nodes, uniformly in $n$), we see that the maps in \eqref{sc-m}
provide uniformly accurate approximations of the maps in \eqref{sc} on
compacta, for $M$ and $n$ sufficiently large.

We consider the quantity
$F^M_{\ogup}(b_{k,j})-F^M_{\ogup}(a_{k,j})=F^M_{\ogup}(a_{k,j}+\ve+\eta_{k,j})-F^M_{\ogup}(a_{k,j})$
(for $M\gg k$) as $n\to\infty$, or equivalently, as $\ve\to0$. In the
integrand, we can introduce the substitution
$z=a_{k,j}+t(b_{k,j}-a_{k,j})$ so that the factor
$$
\f{z-a_{k,j}}{z-b_{k,j}}
$$
becomes $\f{-t}{1-t}$. Moreover, using from the previous passage the
estimate that $|b_{k',j'}-a_{k',j'}|=\ve+o(n)$, the other factors
$$
\f{z-a_{k',j'}}{z-b_{k',j'}}
$$
become $1+\f\ve{c_{k,j}-c_{k',j'}}+o(\ve)$ for $n$ sufficiently large: here
recall that $c_k$ represents the position of the node which is the
limit of the points $a_{k,j}$ and $b_{k,j}$. As a consequence, we compute
that
$$
F^M_{\ogup}(b_{k,j}) - F^M_{\ogdn}(a_{k,j}) = e^{i\pi/4}\int^1_0
\(\f t{1-t}\)^{1/2}dt\(1 + \f\ve2\sum_{\substack{(k',j)\neq(k,j)\\
|k|\le M}}\f1{c_{k,j}-c_{k',j}}\),
$$
where the sum runs over the terms in the integrand of $F^M_{\ogup}$
which are not indexed by $(k,j)$, with $k\le M$.

Now, the computation of $F^M_{\ogdn}(b_{k,j})-F^M_{\ogdn}(b_{k,j})$ is
analogous, yielding
$$
F^M_{\ogdn}(b_{k,j}) - F^M_{\ogdn}(a_{k,j}) = e^{-i\pi/4}\int^1_0
\(\f{1-t}t\)^{1/2}dt\(1 - \f\ve2\sum_{\substack{(k',j)\neq(k,j)\\
|k|\le M}}\f1{c_{k,j}-c_{k',j'}}\).
$$
Since $\int^1_0\(\f{t-1}t\)^{1/2}=\int^1_0\(\f t{1-t}\)^{1/2}dt$, we see
that the horizontal period condition $F_{\ogup}(b_{k,j})-
F_{\ogup}(a_{k,j})=\ov{F_{\ogdn}(b_{k,j})-F_{\ogdn}(a_{k,j})}$
provides that
\begin{equation}\label{forces}
\SF^M_{k,j} = \sum_{\substack{k',j'\neq(k,j)\\ |k|\le M}}
\f1{c_{k,j}-c_{k,j'}} = 0.
\end{equation}
(Evidently, these ``forces'' $\SF^M_{k,j}$ converge to
\begin{equation}\label{force limit}
\SF_{k,j} = \sum^\infty_{M=0}\sum_{\substack{(k',j')\neq(k,j)\\ |k|=M}}
\f1{c_{k,j}-c_{k',j'}} = 0
\end{equation}
but we prefer to continue to work with the approximations for a few more
paragraphs, in order to interpret $\SF_{k,j}$ as a gradient.)
\end{rem}

We have three final goals. We first aim to show that there is a unique
configuration $\{c_{k,j}\}$ which satisfies \eqref{forces}, that this
configuration consists of equally spaced points (to order $o(1)$ in
$M$), and finally that this symmetric configuration is a non-degenerate
critical point.

To begin, we observe that the ``force'' equations \eqref{forces} may be
interpreted as the vanishing of the gradient for the function $\wh
H(\{c_{k,j}\})=-\prod_{\substack{(k,j)\neq(k',j)\\
|k|=M}}|c_{k,j}-c_{k',j'}|^{-1}$. As the $c_{k,j}$ are periodic in $n$ (at
least up to the cut off $M$), we may regard this function as arising
from a function $H$ with domain the simplex
$D=\{0=c_{0,0}<c_{0,1}<\dots<c_{0{\vk}}<1\}$. Clearly the function $H$ is
proper on this simplex, and so we may expect an interior critical point
at an interior global minimum. We now compute the Hessian of $H$, as
equivalently, the differential of the map $\SF^M:D\to\BR^{\vk+1}$ given by
$\SF^M=(\SF^M_{0,0},\dots,\SF^M_{0{\vk}})$. This Hessian has the form
\begin{gather}
\begin{pmatrix}
-\sum\lm_{\substack{(k,j)\neq(0,0)\\ |k|\le M}}(c_{0,0}-c_{k,j})^{-2}
&\sum\lm_{|k|<M}(c_{0,0}-c_{k,1})^{-2}  &\dots&\sum\lm_{|k|\le
M}(c_{0,0}-c_{k,{\vk}})^{-2}\\
\sum\lm_{|k|\le M}(c_{0,1}-c_{k,0})^{-2}  
&-\sum\lm_{\substack{(k,j)\neq(0,1)\\
|k|\le M}}(c_{0,1}-c_{k,j})^{-2}&\dots&\sum\lm_{|k|\le
M}(c_{0,1}-c_{k,{\vk}})^{-2}\\
&&&-\sum\lm_{\substack{(k,j)\neq(0,g)\\ |k|\le M}}(c_{0{\vk}}-c_{k,j})^{-2}
\end{pmatrix}\notag\\
= -\[\sum^{\vk}_{j=0}\sum\lm_{\substack{(k,j)\neq(0,j)\\ |k|\le M}}
(c_{0,j}-c_{k,j})^{-2}\] I - \(\sum\lm_{|k|\le M}
(c_{0,i}-c_{k,j})^{-2}\)^{\vk}_{i,j=0}. \label{hess}
\end{gather}
As each row of the second matrix in \eqref{hess} sums to the negative of
the diagonal entry of the first matrix, we easily see that this Hessian
is negative semi-definite with kernel coming only from a vector
$(\la,\dots,\la)$ with identical entries. As this vector 
$(\la,\dots,\la)$ represents
only a constant translation of the nodes to one direction or other, 
it is not a permissable deformation in $D$, since we required
$c_{0,0}=0$. Thus this Hessian is negative definite on the (projectivized)
domain $D$ of configurations of nodes. Thus each critical point of $-H$
has index $0$, and so, Morse theory applied to the cell $D$ implies that
there is a unique critical point.

Finally let $c^*=(0,\f1{\vk+1},\f2{\vk+1},\dots,\f{\vk}{\vk+1})\in D$ denote the
configuration of equally spaced points in $D$. We observe that
$\SF^M_i(c^*)=o(\f1M)$. Thus, since the unique zero (say $c_M$) of
$\SF^M$ is a minimum of $H$, it is then uniformly bounded away from $\p D$,
and so we see that $c_M\to c^*$ as $M\to\infty$.

We conclude that the equally spaced point set $c^*$ is the unique limit
configuration of the vertices $\{a_{k,j},b_{k,j}\}$. Moreover, from the
analysis above of the Hessians, we observe that $c^*$ is a
non-degenerate zero of the map $\SF_{0,j}:D\to\BR^{\vk+1}$ (where we have
now passed to the limit force equations \eqref{force limit}). Thus, by
the implicit function theorem, there is a unique extension of this
configuration to the unnoded case \eqref{sc}, yielding a zigzag, whose
corresponding minimal surface has Scherk ends making a small positive
angle between them. But as the standard Scherk examples are also such
a family, and the family produced by the implicit function theorem is
unique, we conclude that $\ov{M}_n=\ov{S}(t_n)$ for $n$
large.  This completes the proof of Proposition \ref{small}. \qed

\section{The proofs of Theorems \ref{thm1} and \ref{thm2}.}

We are now in a position to prove Theorem \ref{thm2}.  By the openness
result in
Proposition \ref{diffeo}, the components of $\mathcal{M}_{\vk}$ are curves.
By properness result in Proposition \ref{proper}, for each component $C$
of $\mathcal{M}_{\vk}$, the map 
$\theta \colon C \rightarrow (0, \frac{\pi}{2}]$ is
a diffeomorphism.  By the uniqueness result for small angle in Proposition
\ref{small}, the only component of $\mathcal{M}_{\vk}$ for which the
restriction of $\theta$ is onto is the component of Scherk examples.
Hence, $\mathcal{M}_{\vk}$ consists of the component of Scherk examples,
which proves Theorem \ref{thm2}.

Assume now that $M$ is a connected minimal surface with $A(M)<3\pi$.  In
this case, the limit tangent cone $C$ of $M$ is a cone over an
integral varifold $F$ of $S^2$, consisting of a finite number of geodesic
segments joined at the finite number of vertices of $F$, and at each
vertex
$x_0 \in F$, in a small neghborhood of $x_0$, the varifold $F$ consists of
two geodesics
crossing transversely.  This fact follows immediately from our area growth
assumption and the fact that when considered to be a current, 
the varifold $F$ has no
boundary. From this local property at the vertices, we see immediately
that $F$ is a finite union of circles and our area assumption implies that 
there
are at most two such circles, counted with possible multiplicity.  In
particular, the area growth constant $A(M) = k \pi$, where $k=1$ or $2$.

If $A(M)= \pi$, then $M$ is a plane by the monotonicity formula for
area.  So, assume now that $A(M)=2\pi$. In this case, any limit
tangent cone at infinity for $M$ consists of two planes or a single
plane of multiplicity two.  Now assume that $M$ has
infinite symmetry group and we will prove that $M$ is a catenoid or
one of the Scherk examples.

Since $M$ has infinite symmetry group, then $M$ is invariant under a
continuous group of rotations or it is invariant under a screw motion symmetry.
Assume that $M$ is not a catenoid, which is the only minimal surface
of revolution.  Since $M$ is invariant under a screw motion symmetry,
one sees that it has a unique limit tangent cone $C(M)$ at infinity.
It follows that $C(M)$ consists of two distinct planes or a single
plane of multiplicity two. Since the screw motion symmetry of $M$
leaves $C(M)$ invariant, the screw motion composed with itself four
times is a pure translation $\gamma \colon \rth \rightarrow
\rth$, which, after a rigid motion and homothety, can be taken to be
translation by the vector $(0,0,1)$ that lies in ``both'' planes in
$C(M)$.

Consider the translational subgroup  $\mathbb{Z}=\{(0, 0, n) \mid n
\in \mathbb{N}\}$ of $\rth$, which leaves $M$ invariant.  Note that if
 $\ov{M}=M/\mathbb{Z} \subset \rth/\mathbb{Z}$ has finite
topology,  then Theorem \ref{thm2} implies $M$ is a Scherk
surface.  So, it remains to show $\ov{M}$ has finite topology.  If $C(M)$
consists of two distinct planes $P_1, P_2$, then the facts that each
has multiplicity 1 and $M$ is periodic imply that outside of some
solid cylinder $\Delta$ around $P_1 \cap P_2$, $M$ consists of four graphs
of small gradient over $(P_1 \cup P_2)-\Delta$, which implies $\ov{M}$
has finite topology.

Assume now that $C(M)$ is the $(x_1, x_3)$-plane $P$ (with
multiplicity two as a locally finite integral varifold) and we will
obtain a contradiction.  Since $C(M)$ is unique and $M$ is invariant
under $\mathbb{Z}$, there exists a $\mathbb{Z}$-invariant regular
neighborhood $N(P)$ whose width around $P$ is a positive function
$W(|x_1|)$ which grows sublinearly in the variable $|x_1|$, such that
$M$ is contained in the interior of $N(P)$.  Let $N(A_{x_1})$ in
$\rth/\mathbb{Z}$ be the quotient regular neighborhood of $A_{x_1}$.
Note that every annulus in $\rth/\mathbb{Z}$ which is parallel to
$A_{x_1}$ must intersect $\ov{M}$, otherwise $M$ is contained in a
halfspace in $\rth$ which contradicts the Half Space Theorem
\cite{hm10}, since $M$ is not a plane.  For some $R>0$ large, the
circle $\alpha_R=(0,R) \times \re/\mathbb{Z}$ lies outside of
$N(A_{x_1})$.  For $\theta \in [0, \frac{\pi}{2}]$, let $A(\theta)$ be
the vertical infinite flat half annulus in $[0, \infty) \times [0,
\infty) \times \re/\mathbb{Z}$ with boundary $\alpha_R$ and parallel
to the vector $(\cos(\theta), \sin (\theta), 0)$.

  Without loss of
generality, we may assume that
$A(0) \cap \ov{M} \neq \emptyset$, and so, by the maximum principle
for minimal surfaces, $A(0)$ intersects $\ov{M}$ transversely at some point.  It follows that $A(\theta_0)$
also intersects $\ov{M}$ transversely at some point for some positive
$\theta_0$ close to $0$.  Since $A(\theta_0)$ intersects $N(A_{x_1})$
in a compact set, then there is a nonempty compact subdomain $\Delta$
of $\ov{M}$ which lies in one of the bounded components of
$N(A_{x_1}) - A(\theta_0)$.  Hence, there exists a largest
$\theta_1$ such that $A(\theta_1) \cap \ov{M} \neq \emptyset$ and such
that at every point of this intersection, $\ov{M}$ locally lies on one
side of $A(\theta_1)$.  This contradicts the maximum principle for
minimal surfaces, which completes the proof of Theorem \ref{thm1}.

\bibliographystyle{plain}
\bibliography{bill}
\end{document}